\documentclass{article}[11pt]

\usepackage{tikz, caption, subcaption}
\usepackage[ruled,vlined,linesnumbered]{algorithm2e}
\usepackage{amsfonts,epsf,amsmath,amssymb}
\usepackage{fullpage}

\usepackage{enumerate}
\usepackage{setspace}
\usepackage{color}

\usetikzlibrary{decorations.pathreplacing}

\newtheorem{thm}{Theorem}[section]
\newtheorem{prop}[thm]{Proposition}

\newtheorem{cor}[thm]{Corollary}

\newtheorem{lemma}[thm]{Lemma}

\newcommand{\qed}{\hfill ~$\square$\bigskip}
\newcommand{\proof}{\noindent\textbf{Proof. }}

\def\NN{\hbox{\sf I\kern-.13em\hbox{N}}}
\def\RR{\hbox{\sf I\kern-.14em\hbox{R}}}
\def\ZZ{\hbox{\sf I\kern-.14em\hbox{Z}}}

\newcommand{\vertex}{\node[vertex]}
\tikzstyle{vertex}=[circle, draw, inner sep=0pt, minimum size=6pt]

\newcommand{\gt}{\gamma_t}

\newcommand{\grt}{\gamma_{\rm gr}^t}

\begin{document}

\title{Total dominating sequences in trees, split graphs, and under modular decomposition}

\author{
$^{a,b}$Bo\v{s}tjan Bre\v{s}ar
\and
$^{b}$Tim Kos
\and
$^{c}$Graciela Nasini
\and
$^{c}$Pablo Torres \\
}

\date{\today}

\maketitle

\begin{center}
$^a$ Faculty of Natural Sciences and Mathematics, University of Maribor, Slovenia\\
$^b$ Institute of Mathematics, Physics and Mechanics, Ljubljana, Slovenia\\
$^c$ Depto. de Matem\' atica, Universidad Nacional de Rosario and Consejo Nacional de Investigaciones Cient\' ificas y T\' ecnicas, Argentina\\
\end{center}

\begin{abstract}
A sequence of vertices in a graph $G$ with no isolated vertices is called a total dominating sequence if every vertex in the sequence totally dominates at least one vertex that was not totally dominated by preceding vertices in the sequence, and, at the end all vertices of $G$ are totally dominated (by definition a vertex totally dominates its neighbors). The maximum length of a total  dominating sequence is called the Grundy total domination number, $\grt(G)$, of $G$, as introduced in [B.~Bre{\v{s}}ar, M.~A.~Henning, and D.~F.~Rall, Total dominating sequences in graphs, \textit{Discrete Math.}  \textbf{339} (2016), 1165--1676]. In this paper we continue the investigation of this concept, mainly from the algorithmic point of view. While it was known that the decision version of the problem is NP-complete in bipartite graphs, we show that this is also true if we restrict to split graphs. A linear time algorithm for determining the Grundy total domination number of an arbitrary tree $T$ is presented,  based on the formula $\grt(T)=2\tau(T)$, where $\tau(T)$ is the vertex cover number of $T$. A similar efficient algorithm is presented for bipartite distance-hereditary graphs. Using the modular decomposition of a graph, we present a frame for obtaining polynomial algorithms for this problem in classes of graphs having relatively simple modular subgraphs. In particular, a linear algorithm for determining the Grundy total domination number of $P_4$-tidy graphs is presented. In addition, we prove a realization result by exhibiting a family of graphs $G_k$ such that $\grt(G_k)=k$, for any $k\in{\mathbb{Z}^+}\setminus\{1,3\}$, and showing that there are no graphs $G$ with $\grt(G)\in \{1,3\}$. We also present such a family, which has minimum possible order and size among all graphs with Grundy total domination number equal to $k$.

\end{abstract}




\section{Introduction}

The \emph{total domination number}, $\gt(G)$, of a graph $G$ with no isolated vertices is the smallest cardinality of a set of vertices $S$ such that every vertex of $G$ has a neighbor in $S$. (If the condition
only requires that vertices from $V(G)\setminus S$ have a neighbor in $S$, then the resulting invariant is the \emph{domination number} $\gamma(G)$ of $G$.) 
Let us introduce our main invariant, which is defined for all graphs $G$ without isolated vertices, see~\cite{bhr-2016}.
Let $S=(v_1,\ldots,v_k)$ be a sequence of distinct vertices of $G$. The corresponding set $\{v_1,\ldots,v_k\}$ of vertices from the sequence $S$ will be denoted by $\widehat{S}$. The sequence $S$  is a {\em legal (open neighborhood) sequence} if
\begin{equation}
\label{e:total}
N(v_i) \setminus \bigcup_{j=1}^{i-1}N(v_j) \ne\emptyset.
\end{equation}
holds for every $i\in\{2,\ldots,k\}$.
If, in addition, $\widehat{S}$ is a total dominating set of $G$, then we call $S$ a \emph{total dominating sequence} of $G$. The maximum length of a total dominating sequence in $G$ is called the \emph{Grundy total domination number} of $G$ and denoted by $\grt(G)$; the corresponding sequence is called a \emph{Grundy total dominating sequence} of $G$. 

A motivation for introducing total dominating sequences came from the so-called total domination game~\cite{heklra-2013,heklra-2015}, in which the sequences are a result of two-player game with players having the opposite goals; one player wants the graph to be totally dominated in as few moves as possible, while the other player wants to maximize the sequence of moves. The length of the resulting sequence played in such a game is thus a lower bound for the Grundy total domination number of a graph. A similar game, called the domination game, was introduced earlier with respect to the standard domination number\cite{brklra-2010}, and was already studied in a number of papers.
In particular, motivated by the domination game, a different version of dominating sequences was defined in~\cite{bgm-2014}, in which legality is considered with respect to closed neighborhoods (i.e., in the above definition just replace open neighborhoods by closed neighborhoods in~\eqref{e:total}); longest sequences in that sense are called the {\em Grundy dominating sequences}, and the corresponding invariant the {\em Grundy domination number of a graph}.

Efficient algorithms for the Grundy domination number of trees, cographs and split graphs have been presented in~\cite{bgm-2014}. In addition, minimal dominating sets have been characterized through some algebraic properties of dominating sequences, and some general lower bounds for this parameter were also established. Similarly, a lower bound was obtained for the Grundy total domination number of an arbitrary graph in~\cite{bhr-2016}, with an improvement for $k$-regular graphs. Using the connection with covering sequences in hypergraphs, NP-completeness of the decision version of the Grundy total domination number in bipartite graphs was also established. Nevertheless, no other algorithmic issues were considered in~\cite{bhr-2016}.

In this paper we continue the study of the Grundy total domination number with an emphasis on algorithmic issues. In Section~\ref{sec:bounds} we start setting new bounds for the Grundy total domination number and proving a realization theorem about it. We continue in Section~\ref{sec:modular}, in which we first note, how the two operations on which modular decomposition of a graph depends, namely the join and the disjoint union of two graphs, effect the Grundy total domination number. We follow with an application of these observations in Section~\ref{sec:alg}, by presenting a linear algorithm to determine the Grundy total domination number of $P_4$-tidy graphs. This graph class generalizes several few $P_4$'s graph families, among them cographs, $P_4$-sparse, $P_4$-extendible and $P_4$-reducible graphs \cite{Giako}. 
Besides, we note a first difference between the (standard) Grundy domination number and the Grundy total domination number in their behavior in the class of trees. For the former no explicit formula was found~\cite{bgm-2014}, and the algorithm for determining the Grundy domination number of a tree is based on a recursive, dynamic programming approach. 
In this paper (see Section~\ref{sec:alg}) we prove that the Grundy total domination number of an arbitrary tree can also be found in linear time, but the algorithm is much simpler. It is based on the formula, which expresses this number as two times the vertex cover number of a tree. As it turns out, a similar approach (though no such nice connection with the vertex cover exists) can be used in determining the Grundy total domination number of an arbitrary bipartite distance-hereditary graph, which is done in the second part of the same section.
To conclude, recall that the (standard) Grundy domination number problem was proven to be NP-complete in chordal graphs~\cite{bgm-2014}. On the other hand, in the same paper an efficient algorithm for determining the Grundy domination number of an arbitrary split graphs was presented. To end this paper (see Section~\ref{sec:NP}) we prove that the total version of this problem is NP-complete in split graphs, which shows that the two problems have some essential differences, in spite of the similarity of their definitions. 

\section{Preliminaries}
\label{sec:prelim}

This section is devoted to notational and preliminary issues. For notation and graph theory terminology, we in general follow~\cite{MHAYbookTD}. A \emph{non-trivial graph} is a graph on at least two vertices. 
For each positive integer $n$, $K_n$ and $P_n$ are, respectively, the complete graph and the path with $n$ vertices. For $n\geq 3$, $C_n$ denotes the cycle with $n$ vertices.
Given a graph $G$, $\overline G$ denotes its complement. The subgraph induced by a set $S$ of vertices of $G$ is denoted by $G[S]$. We write $G\setminus S$ for the subgraph induced by $V(G)\setminus S$. If $S=\{v\}$ we simply write $G-v$. 

The \emph{degree} of a vertex $v$ in $G$, denoted $deg_G(v)$, is the number of neighbors, $|N_G(v)|$, of $v$ in $G$. The minimum and maximum degree among all the vertices of $G$ are denoted by $\delta(G)$ and $\Delta(G)$, respectively. 
Given $U\subset V(G)$, $N(U)=\cup_{v\in U} N(v)$ and, for each $k\in \{1,\ldots, |V(G)|\}$, $\delta_k(G)=\min \{|N(U)|\,:\, U\subset V(G), |U|=k\}$. Observe that $\delta_1(G)=\delta(G)$.

A \emph{leaf} is a vertex of degree~$1$, while its neighbor is a \emph{support vertex}. Leaves are also called \emph{pendant vertices}. A \emph{strong support vertex} is a vertex with at least two leaf-neighbors. 
We denote by $L(G)$ the set of leaves of $G$ and for each $v\in V(G)$, $L(v)=L(G)\cap N(v)$. Given $\ell \in L(G)$, $s(\ell)$ is the support vertex of $\ell$, that is, $N(\ell)=\{s(\ell)\}$.

Two vertices $u$ and $v$ in $G$ are called {\em true (resp. false) twin vertices} if $N[u]=N[v]$ (resp. $N(u)=N(v)$). A vertex $v\in V(G)$ is called a {\emph{true (resp. false) twin vertex}} if there exists $u\in V(G)\setminus\{v\}$, such that $u$ and $v$ are true (resp. false) twin vertices.

Given two graphs $G$ and $R$, and $v\in V(G)$, the graph obtained by \emph{replacing} $v$ \emph{by} $R$ \emph{in} $G$ is the graph whose vertex set is  $(V(G)\setminus\{v\})\cup V(R)$ and whose edges either belong to $E(G -v)\cup E(R)$ or connect any vertex in $V(R)$ with any vertex in $N_G(v)$.

If $S=(v_1,\ldots,v_k)$ is a legal (open neighborhood) sequence in a graph $G$ (recall the condition~\eqref{e:total}), then we say that $v_i$ \emph{footprints} the vertices from $N(v_i) \setminus \cup_{j=1}^{i-1}N(v_j)$, and that $v_i$ is the \emph{footprinter} of every vertex $u\in N(v_i) \setminus \cup_{j=1}^{i-1}N(v_j)$. That is, $v_i$ footprints vertex $u$ if $v_i$ totally dominates $u$, and $u$ is not totally dominated by any of the vertices that precede $v_i$ in the sequence. Thus the function $f_S \colon V(G)\to \widehat{S}$ that maps each vertex to its footprinter is well defined. 
Clearly the length $k$ of a total dominating sequence $S$ is bounded from below by the total domination number, $\gamma_t(G)$, of $G$. In particular, $\grt(G)\ge \gamma_t(G)$.

Clearly, a graph with isolated vertices has no total dominating sequences. In order to obtain consistency when working with vertex induced subgraphs, we extend the definition of the Grundy total domination number for such graphs as the maximum length of a legal sequence and a Grundy total dominating sequence 
is a legal sequence of maximum length. 
Clearly, for graphs without isolated vertices this definition coincides with the one given in ~\cite{bhr-2016}. Moreover, if $E(G)=\emptyset$, there are no legal sequences except the \emph{empty sequence}, denoted by $S=()$ and in this case, we let $\grt(G)=0$. We also define the parameter $\eta(G)$ having value zero if $G$ has no isolated vertices and one, otherwise. This notion, as well as the extension of the Grundy total domination number to disconnected graphs will be used in Section~\ref{sec:modular}.

Let $S_1=(v_1,\ldots , v_n)$ and $S_2=(u_1,\ldots , u_m)$, $n,m \geq 0,$ be two sequences in $G$, with $\widehat S_1\cap \widehat S_2=\emptyset$. The {\em concatenation} of $S_1$ and $S_2$ is defined as the sequence $S_1 \oplus S_2=(v_1,\ldots , v_n,u_1,\ldots , u_m)$. Clearly $\oplus$ is an associative operation on the set of all sequences, but is not commutative. Moreover, $()\oplus S= S \oplus ()= S$, for any sequence $S$.

A legal sequence $S$ in a graph is said to be \emph{maximal} if $S\oplus S'$ is a legal sequence only if $S'=()$.
Clearly, every Grundy total dominating sequence is maximal and thus, the Grundy total domination number is the maximum length of a maximal legal sequence.


\section{New bounds and realization results for $\grt(G)$}\label{sec:bounds}

\subsection{Bounds}

For a matching $M$ in a graph $G$ a vertex incident to an edge of $M$ is called \emph{strong} if its degree is~$1$ in the subgraph $G[V(M)]$.  The matching $M$ is called a \emph{strong matching} (also called an \emph{induced matching} in the literature) if every vertex in $V(M)$ is strong. The number of edges in a maximum (strong) matching of $G$ is the \emph{(strong) matching number}, ($\nu_s(G)$) $\nu(G)$, of $G$. The strong matching number is studied, for example, in~\cite{JRS-2014,KMM-2012}. As defined in~\cite{GH-2005}, $M$ is a \emph{semistrong} matching if every edge in $M$ has a strong vertex. The number of edges in a maximum semistrong matching of $G$ is the \emph{semistrong matching number}, $\nu_{ss}(G)$, of $G$.

We recall the lower bound of the Grundy total domination number based on the semistrong matching number:

\begin{prop}
{\rm \cite{bhr-2016}}
For every graph $G$, $\grt(G) \ge 2\nu_{ss}(G)$.
\label{indmatch}
\end{prop}

Next, we present an upper bound based on the vertex cover number for arbitrary graphs.
Given a graph $G$, let $\tau(G)$ denote the {\em vertex cover number}, 
that is, the minimum cardinality of a set $S$ of vertices in $G$ such that each edge in $G$ is incident to a vertex from $S$.

\begin{prop}
For every graph $G$, $\grt(G) \le 2\tau(G)$.
\label{prp:tau}
\end{prop}
\proof
Let $S$ be a Grundy total dominating sequence of $G$,
and let $u$ be an arbitrary vertex from a maximum vertex cover $C$ of $G$. 
Clearly, exactly one of the vertices from $N(u)\cap \widehat{S}$ footprints 
$u$. Also, $u$ is a footprinter of some vertices from $N(u)$ if and only if $u\in \widehat{S}$. In other words, for each $u\in C$, the edges incident with $u$ are involved in footprinting operation at most twice, from which the desired bound follows. 
\qed

By K\"{o}nig's theorem, in bipartite graphs the minimum cardinality of a vertex cover coincides with that of a maximum matching. Thus, using also Proposition~\ref{indmatch}, we get the following bound.

\begin{cor}\label{c:bip}
For every bipartite graph $G$, $2\nu_{ss}(G) \le \grt(G) \le 2\nu(G)$.  
\end{cor}

The following result shows how the parameters $\delta_k(G)$ are used to obtain new upper bounds for the Grundy total domination number.
\begin{lemma}\label{l:delta}
For each $k\in \{1,\ldots, |V(G)|\}$, $\grt(G)\leq k+|V(G)|-\delta_k(G)$.
\end{lemma}
\proof
If $k\geq\grt(G)$ the result immediately follows. Let $k<\grt(G)$ and $S=(x_1,\ldots,x_k,\ldots,x_t)$ a Grundy total dominating sequence of $G$. Consider $U=\{x_1,\ldots,x_k\}$. Clearly, $\grt(G)-k=t-k\leq  |V(G)|-|N(U)|\leq |V(G)|-\delta_k(G)$. Then, $\grt(G)\leq k+|V(G)|-\delta_k(G)$.
\qed

Note that $\grt(K_n)=1+n-\delta_1(G)=2+n-\delta_2(G)=2$. Moreover, for $P_n$, with $n$ even, and $C_n$ with $n$ odd, the Grundy total domination number coincides with these upper bounds, for all $k\in\{1,\dots,n\}$.
On the other hand, observe that there exist graphs with Grundy total domination number strictly smaller than these upper bound for every $k\in\{1,\dots,|V(G)|\}$. For instance, consider $P_n$ with $n$ odd. In these cases $\delta_k(P_n)\leq k$ for every $k\in\{1,\dots,n\}$ but $\grt(P_n)=n-1$.

\medskip

Next, the following bounds concern the deletion of an arbitrary vertex of a graph.

\begin{lemma}\label{grundysubgraph}
If $G$ is a graph, $v\in V(G)$ and $G'=G-v$, then $$\gamma_{gr}^t(G)-2\leq\gamma_{gr}^t(G')\leq\gamma_{gr}^t(G).$$
\end{lemma}
\proof
The upper bound follows from the fact that every legal sequence of $G-v$ is a legal sequence of $G$.
On the other hand, if $S$ is a total dominating sequence of $G$ and $u$ is the vertex that footprints $v$, then $S\setminus\{u,v\}$ is a total dominating sequence of $G-v$ and the lower bound follows. 
\qed

To see that the lower bound is tight note that if $G$ is the net graph (i.e., the graph obtained from the triangle by adding a pendant vertex to each vertex of the triangle), then $\gamma_{gr}^t(G)=\gamma_{gr}^t(G-v)+2=6$ where $v$ is a pendant vertex of $G$. It is also easy to find examples where the upper bound in Lemma~\ref{grundysubgraph} is attained. In particular, considering false twin vertices we have the following result.

\begin{prop}
\label{cor:twins}
Let $G$ be a graph and $v,v'$ a pair of false twins in $G$. Let $G'=G-v'$ and $S'$ be a Grundy total dominating sequence in $G'$.
Then, $S'$ is a Grundy total dominating sequence in $G$ and $\gamma_{gr}^t(G)=\gamma_{gr}^t\left(G-v'\right)$.
\end{prop}
\proof
The results follow from the facts that for every legal sequence $S$, we have $|\widehat S\cap\{v,v'\}|\leq1$, and that if $w\in \widehat S$ footprints $v$, then $w$ also footprints $v'$.
\qed

We finally study the total domination number of graphs with pendant vertices.
Note that if $v$ is a vertex in $G$ such that $L(v)\neq\emptyset$, then $v$ belongs to every maximal legal sequence of $G$. Moreover, the following result can be proven for vertices with leaf-neighbors. 

\begin{prop}
\label{cor:pendant}
Let $G$ be a graph, $v\in V(G)$ such that $L(v)\neq\emptyset$, $G'=G\setminus(\{v\}\cup L(v))$, and $S'$ a Grundy total dominating sequence of $G'$. Then, for any $\ell\in L(v)$, $S= (\ell)\oplus S'\oplus (v)$ is a Grundy total dominating sequence of $G$ and  $\gamma_{gr}^t(G)=\gamma_{gr}^t(G')+2$.
\end{prop}

\proof
Let $G,G'$ and $S'$ be defined as in the statement of the proposition. Clearly, $(l)\oplus S'\oplus (v)$ is a legal sequence of $G$ for every $l\in L(v)$. Hence, $\gamma_{gr}^t(G)\geq \gamma_{gr}^t(G')+2$.
On the other hand, if $T$ is an arbitrary maximal legal sequence of $G$ and $u$ the footprinter of $v$ in $T$, then $T\setminus\{u,v\}$ is a legal 
sequence of $G'$. Thus, $|\widehat{T}|\leq \gamma_{gr}^t(G')+2$, which implies $\gamma_{gr}^t(G)\leq \gamma_{gr}^t(G')+2$. Since $|\widehat S|= \gamma_{gr}^t(G')+2$, $S$ is a Grundy total dominating sequence of $G$.
\qed

\subsection{Realization of Grundy total domination numbers}
\label{sec:realization}

Since for any graph $G$ we have $2\le \gamma_t(G)\le \grt(G)$, it is clear that there are no graphs $G$ with Grundy total domination number equal to $1$. We will prove that there are also no graphs $G$ such that $\grt(G)=3$, and one may wonder, whether there are any graphs $G$ such that their Grundy total domination number is odd. (It is easy to find graphs with even total domination numbers of arbitrary size, for instance among paths, because $\grt(P_n)=2\lfloor n/2\rfloor$.) As we will show next, number $3$ is exceptional here, since all other odd numbers greater than 3 can be realized as the total domination numbers of some graphs. 

First, let us present a result, which shows, why even numbers are in some sense much easier to be realized as total domination numbers of graphs than the odd numbers. 

\begin{prop}
\label{maxind}
Let $G$ be a graph, and $S$ a Grundy total dominating sequence. Suppose that $S=S_1\oplus S_2$, and $\widehat S_1$ is a maximal independent set in $G$.
Then $|\widehat S_2|=|\widehat S_1|$, and so $\grt(G)$ is even.
\end{prop}
\proof
Let $S=(s_1,\ldots,s_r,s_{r+1},\ldots,s_t)$, where $\widehat S_1=\{s_1,\ldots,s_r\}$ is a maximal independent set in $G$. Note that $S_1$ totally dominates exactly the vertices from $V(G)\setminus \widehat S_1$, and so $|\widehat S_2|=t-r$ is at most $r$. We claim that $|\widehat S_2|=r$.

Note that a vertex $u_i$ from $f^{-1}_S(s_i)$ (i.e., a vertex footprinted by $s_i$ within the sequence $S$) is not adjacent to any of the vertices from $\{s_1,\ldots,s_{i-1}\}$, where $i\in \{1,\ldots, r\}$. Hence the sequence $(s_1,\ldots,s_r,u_r,u_{r-1},\ldots,u_1)$ is a legal open neighborhood sequence of the largest possible length. Since $S$ is a Grundy dominating sequence of $G$, we derive that $|S_2|=r$. 
\qed

Now, we are ready to prove that there are no graphs $G$ with Grundy total domination number equal to $3$.

\begin{prop}
\label{grt3}
There exists no graph $G$ such that $\grt(G)=3$.
\end{prop}
\proof
Suppose that $S=(s_1,s_2,s_3)$ is a Grundy total dominating sequence of a graph $G$.
If $s_1$ and $s_2$ are not adjacent, the corresponding set $\widehat S_1=\{s_1,s_2\}$ is a subset of a maximal independent set $I$ in $G$. We derive, by using Proposition~\ref{maxind} and its proof, that there exists a total dominating sequence in $G$ of legth $2|I|$, which is clearly greater than 3, a contradiction. Hence, we may assume that $s_1$ and $s_2$ are adjacent, and let $u$ be a vertex footprinted by $s_3$. Then $(u,s_1,s_2,s_3)$ is a legal sequence, because $u$ footprints $s_3$, $s_1$ footprints $s_2$, $s_2$ footprints $s_1$ and $s_3$ footprints $u$, a contradiction with $S$ being a Grundy total dominating sequence. We infer that there is no such graph $G$.
\qed

Moreover, we prove that exactly the integers from ${\mathbb{Z}^+}\setminus\{1,3\}$ can be realized as Grundy total domination numbers of some graphs.

\begin{thm} 
For any $n\in {\mathbb{Z}^+}\setminus\{1,3\}$ there exists a graph $G_n$ such that $\grt(G_n)=n$. 
\end{thm}
\proof 
Let $G_n$ be the prism $K_n\Box K_2$, and denote the vertices of one $n$-clique by $a_1,\ldots, a_n$ and of the other $n$-clique by $b_1,\ldots,b_n$ and for each $i$ let $a_i$ be adjacent exactly to $b_i$ among all vertices from the other clique.

Observe that $G_2=C_4$ and so $\grt(G_2)=2$. 

Consider now $G_n$ for $n\ge 4$. It is easy to see that $(a_1,\dots,a_n)$ is a legal sequence of $G_n$. Besides, notice that $\delta_2(G_n)=n+2$. From the upper bound given in Lemma \ref{l:delta} we have that $\grt (G_n)\leq 2+2n-\delta_2(G_n)=n$ and the result follows.
\qed

In addition,  we now consider the problem of finding connected graphs $G_n$ with minimum number of vertices and edges verifying $\grt (G_n)=n$. Clearly, if $n$ is even, $G_n=P_n$.

If $n$ is odd, we know that $n\geq 5$. From \cite[Theorem 4.2]{bhr-2016} it is known that if a graph $G$ verifies $\grt (G)=|V(G)|$ then $|V(G)|$ is even, and so $|V(G_n)|\geq n+1$.
Considering the number of edges, it can be proved that $|E(G_n)|\geq |V(G_n)|+1$. In fact, since $G_n$ is connected and the Grundy total domination number of trees is even (see Theorem~\ref{thm:trees}, we have that $|E(G_n)|\geq |V(G_n)|$.
However, if $|E(G_n)|=|V(G_n)|$, then $G_n$ is a tree with an additional edge defining a cycle. Applying Proposition~\ref{cor:pendant} and knowing that the Grundy total domination numbers of cycles and paths are even, it follows that if $|E(G_n)|=|V(G_n)|$, $\grt (G_n)$ is even. Thus, $|E(G_n)|\geq |V(G_n)|+1$, as claimed.

For $n=5$, we consider the graph $G_5$ with $6$ vertices defined by two disjoint triangles and one edge between one vertex in each triangle.
This graph clearly has the minimum number of vertices and edges.

\begin{figure}[h!]
\begin{center}
\begin{tikzpicture}[]
\tikzstyle{vertex}=[circle, draw, inner sep=0pt, minimum size=6pt]
\tikzset{vertexStyle/.append style={rectangle}}
	\vertex (1) at (-2,1.3) {};
	\vertex (2) at (-2,-1.3) {};
	\vertex (3) at (-0.7,0) {};
	\vertex (4) at (0.7,0) {};
	\vertex (5) at (2,1.3)  {};
	\vertex (6) at (2,-1.3) {};
	\vertex (7) at (3.5,1.3)  {};
	\vertex (8) at (6,1.3)  {};
	\draw [decorate,decoration={brace,mirror, amplitude=10pt,raise=8pt}] (3.5,1.3) -- (6,1.3) node [midway, below=19pt,xshift=1pt] {$P_{2k}$};
	\path
		(1) edge (2)
		(2) edge (3)
		(1) edge (3)
		(3) edge (4)
		(4) edge (5)
		(4) edge (6)
		(5) edge (6)
		(5) edge (7)
		(7) edge[dashed] (8)

	;
\end{tikzpicture}
\end{center}
\caption{Graphs $G_{5+2k}$.}
\label{fig:G_n}
\end{figure}
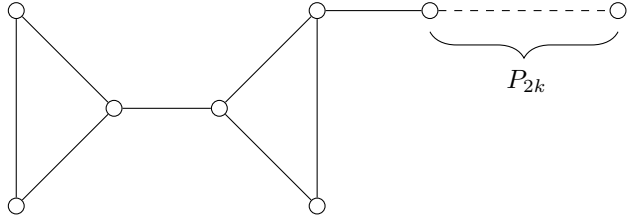

Now, for $n=5+2k$, let $G_n$ be the graph obtained from $G_5$ by adding $P_{2k}$ to a vertex of $G_5$ with degree 2 (see Figure \ref{fig:G_n}).
From Proposition \ref{cor:pendant}, we have that $\grt(G_n)=\grt(G_5)+2k= 5+2k=n$, and so $G_n$ has minimum number of vertices and edges.


\section{Grundy total dominating sequences under modular decomposition}
\label{sec:modular}

Clearly, if a graph $G$ (resp. $\overline G$) is not connected, it can be obtained by disjoint union (resp. join) of two non-empty graphs. If $G$ and its complement are connected, we say that $G$ is \emph{modular}. Given a graph family $\mathcal F$, we denote by $M(\mathcal F)$ the family of modular graphs in $\mathcal F$.

Let us first analyze the behaviour of the Grundy total domination number under disjoint union and join of two graphs. 
Given non-empty graphs $G_1$ and $G_2$ with disjoint sets of vertices, $G_1\oplus G_2$ denotes their disjoint union and $G_1\vee G_2$ their join. (Recall that in this section we consider Grundy total dominating sequence to be the longest legal open neighborhood sequence, which coincides with the definition in graphs without isolated vertices.)

The result for the disjoint union of graphs easily follows:

\begin{prop} \label{union}
Let $G_1$ and $G_2$ be non-empty graphs with disjoint set of vertices.
If $S_1$ and $S_2$ are Grundy total dominating sequences of $G_1$ and $G_2$, respectively, then $S_1\oplus S_2$ is a Grundy total dominating sequence of $G_1\oplus G_2$. Therefore, $\grt(G_1\oplus G_2) = \grt(G_1) + \grt(G_2)$.  
\end{prop}

Considering the join of graphs, we need to distinguish between the graphs with and those without isolated vertices. In this line, remind that the parameter $\eta(G)$ has value one if the graph has isolated vertices and zero, otherwise.

\begin{prop} \label{join}
Let $S_1$ and $S_2$ be, respectively, Grundy total dominating sequences of non-empty graphs $G_1$ and $G_2$ with disjoint set of vertices. Let $G= G_1\vee G_2$.
\begin{enumerate}
	\item If $\eta(G_1)=\eta(G_2)=0$ and $|\widehat S_1|\geq |\widehat S_2|$, $S_1$ is a Grundy total dominating sequence of $G$. 
	\item If $\eta(G_1)=\eta(G_2)=1$, $|\widehat S_1|\geq |\widehat S_2|$ and $v$ is an isolated vertex of $G_1$,  then $S=(v)\oplus S_1\oplus (w)$ is a Grundy total dominating sequence of $G$, for any $w\in V(G_2)$.  
	\item If $\eta(G_1)=1$, $\eta(G_2)=0$, $|\widehat S_1|+2\geq |\widehat S_2|$ (resp. $|\widehat S_1|+2 \leq |\widehat S_2|-1$) then $S=(v)\oplus S_1\oplus (w)$, for any isolated vertex $v$ of $G_1$ and any $w\in V(G_2)$, (resp. $S=S_2$) is a Grundy total dominating sequence of $G$.
\end{enumerate}

Hence, $$\grt (G)=max\{\grt (G_1)+2 \eta(G_1), \grt (G_2)+2 \eta(G_2)\}.$$

\end{prop}

\proof
First observe that, for $i=1,2$, $S_i$ is a legal sequence of $G$. Moreover, if $G_i$ has an isolated vertex $v$ then $(v)\oplus S_i\oplus (w)$ is a legal sequence of $G$, for every  $w\in V(G_j)$, $j=1,2$, $j\neq i$. Therefore,
$$\gamma^t_{gr}(G)\geq max\{\grt(G_1)+ 2 \eta(G_1), \grt(G_2)+ 2 \eta(G_2)\}.$$

In order to prove the three items, it only remains to prove the opposite inequality.

Let $S=(x_1,\ldots,x_{k-1},x_k)$ be a legal sequence of $G$. 
If $\widehat S\cap V(G_i)\neq \emptyset$ for $i=1,2$,  without loss of generality we can assume that $x_1\in V(G_2)$. Let $j=\min\{i: 2\leq j \leq k \;,\; x_i\in V(G_1)\}$. Since $\cup_{i=1}^j N(x_i)=V(G)$, $j=k$. Then, $\widehat S\cap V(G_1)= \{x_k\}$.

If $v$ is an isolated vertex of $G_2$ and $v\in \widehat S$, $v=x_1$ and $(x_2,\ldots,x_{k-1})$ is a legal sequence of $G_2$. Otherwise, if no isolated vertex of $G_2$ belongs to $\widehat S$, $(x_1,\ldots,x_{k-1})$ is a legal sequence of $G_2$. In both cases, $k\leq \grt(G_2)+ 2 \eta(G_2)$. 

It only remains to consider the cases $\widehat S\cap V(G_i)= \emptyset$ for some $i\in\{1,2\}$. Following the same reasoning as before, if $v$ is an isolated vertex of $G_i$ and $v\in \widehat S$, $v=x_1$ and $(x_2,\ldots,x_k)$ is a legal sequence of $G_i$. Otherwise, $(x_1,\ldots,x_{k})$ is a legal sequence of $G_i$. In both cases the inequality 
$k\leq \grt(G_i)+ 2 \eta(G_i)$ is satisfied.
\qed

We introduce the concept of \emph{a modular decomposition tree} of a graph $G$, slightly different to that defined by V. Giakoumakis et al. in \cite{Giako}.

\medskip

A \emph{labeled rooted complete binary tree} is a triple $(T, v, L)$ where $T$ is a complete binary tree, $v$ is an internal vertex of $T$ and $L=\{l(w): w\in V(T)\}$ is a list of labels associated with its vertices. Given $(T_1, v_1, L_1)$ and $(T_2, v_2, L_2)$  with $V(T_1)\cap V(T_2)=\emptyset$ and $v\notin V(T_1)\cup V(T_2)$, we define $(T_1, v_1, L_1) \; v \; (T_2, v_2, L_2)$ as the labeled rooted complete binary tree $(T,v,L)$ such that $V(T)= V(T_1)\cup V(T_2) \cup \{v\}$, $E(T)= E(T_1)\cup E(T_2) \cup \{(v,v_1),(v,v_2)\}$ and $L=L_1\cup L_2\cup \{l(v)\}$.

A \emph{graph decomposition tree} is a labeled rooted complete binary tree such that the labels corresponding with the leaves are vertex disjoint graphs and those corresponding with internal nodes belong to $\{\oplus , \vee\}$. 
Given a graph $G$, a \emph{modular decomposition tree} $T(G)$ \emph{of} $G$ is a graph decomposition tree
constructed recursively as Algorithm \ref{al:MDT} shows.

\begin{algorithm}[hbt!]
\caption{$MDT(G,T,v,L)$}
\label{al:MDT}
\KwIn{$G$, a graph.}
\KwOut{$(T,v,L)$, a modular decomposition tree of $G$.}
\BlankLine
{
\If {$G$ is modular}  {$v$ is an arbitrary vertex, $T=({v},\emptyset)$, $L=\{G\}$\\STOP}
\If {$G$ is not connected} {$*=\oplus$}
\If {$\overline G$ is not connected} {$*=\vee$}
Let $G_1$, $G_2$ be two graphs such that $G=G_1 * G_2$\\
Execute $MDT(G_1,T_1,v_1,L_1)$\\
Execute $MDT(G_2,T_2,v_2,L_2)$ \\
Choose $v\notin V(T_1)\cup V(T_2)$\\
$l(v)=*$\\
$(T,v,L)=(T_1,v_1,L_1) \; v \; (T_2,v_2,L_2)$
}
\end{algorithm}

It is not hard to see that Algorithm \ref{al:MDT} is linear (see \cite{Giako}). Also, the  postorder traversal of $T(G)$ is done in $O(|V(T(G))|)$ time and $|V(T(G))|< 2 |V(G)|$. 

Moreover, from Propositions \ref{union} and \ref{join} we have that Grundy total dominating sequences of $G_1\oplus G_2$ and $G_1\vee G_2$ can be obtained in linear time from Grundy total dominating sequences of $G_1$ and $G_2$. 
Therefore, given a graph class $\mathcal F$, if the problem of finding a Grundy total dominating sequence for graphs in $M(\mathcal F)$ is polynomial (linear) then it is polynomial (linear) for graphs in $\mathcal F$.

In particular, if $\mathcal F$ is the family of cographs, i.e., the graphs that do not contain $P_4$ as an induced subgraph, it is known that $M(\mathcal F)$ only contain trivial graphs for which $S=()$ is the Grundy total dominating sequence. This implies the following result.

\begin{thm}
The Grundy total domination number can be obtained in linear time for cographs.
\end{thm}

Several generalizations of cographs have been defined in the literature, such as $P_4$-sparse \cite{Hoang}, $P_4$-lite \cite{JO3}, $P_4$-extendible \cite{JO4} and $P_4$-reducible graphs \cite{JO5}. A graph class generalizing all of them is the class of $P_4$-tidy graphs \cite{Giako}, and we deal with this class of graphs, by using the mentioned approach, in the subsection \ref{spiders}.

\section{Efficient algorithms for Grundy total dominating sequences}
\label{sec:alg}

\subsection{Trees}

In this section we present an efficient algorithm to determine 
the Grundy total domination number of an arbitrary tree.
As opposed to the algorithm for the (non-total) Grundy domination number of a tree
as presented in~\cite{bgm-2014}, which was very involved, the algorithm we present here is quite straightforward. The algorithm is based on the formula that connects the Grundy total domination number of a tree to its vertex cover number.

\begin{algorithm}[hbt!]
\caption{Grundy total dominating sequence of a tree with no isolated vertices.}
\label{al:grundy}
\KwIn{A tree $T$.}
\KwOut{A Grundy dominating sequence $S$ of $T$.}
\BlankLine
{
	    $S=()$;\\
		$T'$ = T;\\
		$i=0$;\\
	\While{$T'$ has non-isolated vertices}
	{
	    Choose a leaf $u\in T'$, and let $\{w\}=N_{T'}(u)$;\\
		$S=S\oplus (u)$;\\
		$i=i+1$;\\ 
		$w_i=w$;\\
		$T'=T'\setminus\{u,w\}$
		}
	\While{$i>0$}
	{
	$S=S\oplus (w_i)$;\\
	$i=i-1$;\\
		}
}
\end{algorithm}

\begin{thm}
\label{thm:trees}
Algorithm~\ref{al:grundy} returns a Grundy total dominating sequence of an arbitrary tree $T$, with length twice the vertex cover number of $T$. In particular, for any tree $T$, $$\grt(T)=2\tau(T).$$
The complexity of Algorithm~\ref{al:grundy} is $O(|V(T)|)$. 
\end{thm}
\proof
Let us start the proof by showing that the sequence $S$ produced by Algorithm~\ref{al:grundy} is legal; i.e. every vertex in the sequence footprints some vertex in $T$.  Clearly, when a leaf $u$ of a tree $T'$ is chosen in some step of the algorithm, it is clear that its unique neighbor $w$ in $T$ has not yet been totally dominated in previous step (this is because in each step of the first WHILE loop, only one vertex becomes totally dominated and then it is removed from $T'$ together with the leaf that footrprinted it). This implies that $u$ footprints $w$ (for later purposes denoted also by $w_i$), and so the first half of the sequence $S$ is legal. The second half is constructed from vertices $w_i$ that are footprinted in the first part, only that they are listed in the reversed order. At the time vertex $w_i$ appears in $S$, its neighbor $u$ (let us denote it by $u_i$) that footprinted $w_i$, is not yet dominated. Indeed, let $T'$ be the tree at the time rigth before $u_i$ is added to the sequence. Clearly, $u_i$ is a leaf of $T'$, and $w_i$ is its support vertex. Now, all vertices $w_j$ that appear in $S$ before $w_i$ (if any) are from $T'$. This proves that the vertex $w_i$ footprints the vertex $u_i$, that footprinted $w_i$, i.e. $f_S(w_i)=u_i$ and $f_S(u_i)=w_i$. This yields that $S$ is a legal sequence. 

The proof of the correctness and also of the formula $\grt(T)=2\tau(T)$ is based on the 
Proposition~\ref{prp:tau}, which states that twice the vertex cover number is an upper bound for the Grundy total domination number of any graph. We will prove that the set $W$ of vertices $w_i$ (support vertices of the chosen leaves) that are produced by the algorithm forms a vertex cover of the tree $T$. Since the length of the sequence $S$ produced by the algorithm is twice the cardinality of $W$, this implies that $|\widehat{S}|\ge 2\tau(T)$. Now, combining this with Proposition~\ref{prp:tau} we infer that $|\widehat{S}|=2\tau(T)$, which at the same time also implies that $S$ is indeed the longest possible legal sequence, hence $S$ is a Grundy total dominating sequence of $T$. 

To complete the proof we thus need to show that the set $W$ of vertices $w_i$ from $S$ is a vertex cover of $T$. Note that at the end of Algorithm~\ref{al:grundy} the remaining tree $T'$ is either empty or it consists of the set $I$ of isolated vertices. Since $T$ has no isolated vertices, each vertex from $I$ has a neighbor in $T$. We claim that each of the neighbors of a vertex $x$ from $I$ is in $W$. Suppose that $x$ has a neighbor $u$ that is not from $W$. Then, as $u$ is in $S$, it must be in the first half of $S$, i.e., at the time it was added to $S$, $u$ was a leaf of some subtree $T'$ of $T$. But then $u$ had only one neighbor in $T'$, which is a vertex in $W$. This implies that all edges, incident with vertices from $I$ are also incident with a vertex from $W$. Now, if $u$ is a vertex from the first half of $S$, then $u$ is not adjacent to some $u'$, which is also in the first half of $S$. Indeed, $u'$ cannot appear in $S$ after $u$, because then $u'$ would also be in the tree $T'$ at the time $u$ is added to $S$; but at that time, $u$ is adjacent only to a vertex $w$ in $T'$ (where $w\in W$). Now, if $u'$ appears in $S$ before $u$, then again, by arguing in the same way, we find that $u$ is not adjacent to $u'$. Hence, each edge of $T$ is incident to at least one vertex from $W$. This implies that $W$ is a vertex cover of $T$, which completes the proof (time complexity $O(|V(T)|)$ is obvious).
\qed

\subsection{Bipartite distance-hereditary graphs}

A graph $G$ is \emph{distance-hereditary} if for each induced connected subgraph $G'$ of $G$ and all 
$x,y\in V(G')$, $d_{G'}(x,y)=d_G(x,y)$, where $d_{G'}(x,y)$ is the distance in $G'$ between $x$ and $y$, i.e. the length of a shortest path in $G'$ between $x$ and $y$. Then, a graph is \emph{bipartite distance-hereditary} if it is distance-hereditary and bipartite.

It is known that a graph $G$ is distance-hereditary if and only if it can be constructed from 
$K_1$ by a sequence of three operation: adding a pendant vertex, creating a true twin vertex and creating a false twin vertex \cite{bm-1986}.

In \cite{HF-1990} Hammer and Mafray defined a \emph{pruning sequence} of a graph $G$ as a total ordering $\sigma = [x_1,\dots ,x_{|V(G)|}]$
of $V(G)$ and a sequence $Q$ 
of triples  $q_i=(x_i, Z, y_i)$ for $i=1,\ldots, |V(G)|-1$, where $Z\in \{P,F,T\}$ and such that, for $i\in\{1,\dots,|V(G)|-1\}$, if $G_i=G\setminus \{x_1,\dots ,x_{i-1}\}$ then,
\begin{itemize}
 \item $Z=P$, if $x_i$ is a leaf and $y_i=s(x_i)$ in $G_i$,
 \item $Z=F$, if $x_i$ and $y_i$ are false twins in $G_i$,
 \item $Z=T$, if $x_i$ and $y_i$ are true twins in $G_i$.
\end{itemize}

In \cite{HF-1990} distance hereditary graphs are characterized as the graphs which admit a pruning sequence. Later, in \cite{DHP-2001} it is showed that distance hereditary graphs can be recognized in $O(|V(G)|+|E(G)|)$ and given a distance hereditary graph $G$, a pruning sequence of $G$ can be computed in $O(|V(G)|+|E(G)|)$.

Analogously, bipartite distance-hereditary graphs are characterized as the graphs that can be constructed from $K_1$ by a sequence of additions of false twin and pendant vertices. Then, a pruning sequence of a bipartite distance hereditary graph has no words $(x,T,y)$.

From Propositions~\ref{cor:twins} and~\ref{cor:pendant} we derive the recursive Algorithm~\ref{al:grundy2}, which determines a Grundy total dominating sequence of an arbitrary bipartite distance-hereditary graph.

\begin{algorithm}[hbt!]
{\caption{GrundyBDH($G,S)$}
\label{al:grundy2}
\KwIn{A bipartite distance-hereditary graph $G$.}
\KwOut{A Grundy total dominating sequence $S$ of $G$.}
\BlankLine
{
     \If {$E(G)=\emptyset$} {$S=()$ \\ STOP.}
     Obtain a pruning sequence 
     $Q=[q_1,\dots,q_{|V(G)|-1}]$ of $G$\\
     \For {$i=1$ to $|V(G)|-1$}
     {
     \If {$q_i\neq (x_i, F, y_i)$}
      {$G'=G\setminus\{x_1,\ldots,x_i, y_i\}$\\
      Execute GrundyBDH($G'$,$S'$)\\
      $S=(x_i)\oplus S' \oplus (y_i)$\\
      STOP.}
}
}
}
\end{algorithm}

\begin{thm}
Algorithm~\ref{al:grundy2} returns a Grundy total dominating sequence of an arbitrary bipartite distance-hereditary graph $G$. The complexity of Algorithm~\ref{al:grundy2} is $O(|V(G)|(|V(G)|+|E(G)|))$.
\end{thm}
\proof
Note that if $E(G)\neq\emptyset$ there exists $i\in \{1,\ldots, |V(G)|-1\}$ such that $q_i=(x_i, P, y_i)$. Then, Proposition~\ref{cor:pendant} and Proposition~\ref{cor:twins} guarantee the correctness of the algorithm.

Finally, as we have mentioned, a pruning sequence of a (bipartite) distance hereditary graph $G$ can be computed in $O(|V(G)|+|E(G)|)$ \cite{DHP-2001}. Thus, the time complexity $O(|V(G)|(|V(G)|+|E(G)|))$ follows from the fact that step in line 5 runs at most $|V(G)|-1$ times.
\qed

Forest graphs are bipartite distance hereditary graphs for which proposition ~\ref{cor:pendant} allows us to simplify Algorithm~\ref{al:grundy2} (see Algorithm~\ref{al:grundy1}).

\begin{algorithm}[hbt!]
\caption{GrundyForest$(T,S)$
\label{al:grundy1}}
\KwIn{A forest $T$.}
\KwOut{A Grundy total dominating sequence $S$ of $T$.}
\BlankLine
{
     \If {$E(T)=\emptyset$} {$S=()$,\\ STOP.}
     Choose $\ell \in L(T)$ and let $T'=T\setminus ({s(\ell)} \cup L(s(\ell)))$\\
     Execute GrundyForest($T'$,$S'$)\\
     $S=(\ell)\oplus S' \oplus (s(\ell))$\\ STOP.   \\
}
\end{algorithm}

Note that, when the input is a tree, Algorithm~\ref{al:grundy1} uses the same steps, and returns the same sequence of vertices as Algorithm~\ref{al:grundy}. The only difference is that the former algorithm keeps isolated vertices, while the latter deletes them, when they appear as leaves of the support vertex $s(\ell)$.

\subsection{$P_4$-tidy graphs.}\label{spiders}

We start by defining $P_4$-tidy graphs. 
Let $U$ be a subset of vertices inducing a $P_4$ in $G$. A \emph{partner} of $U$ is a vertex $v\in G-U$ such that $U\cup\{v\}$ induces at least two $P_4$s in $G$. 
A graph $G$ is $P_4$-tidy if any $P_4$ has at most one partner. It is known that the class of $P_4$-tidy graphs is self-complementary and hereditary \cite{Giako}.

Non-trivial modular $P_4$-tidy graphs are the graphs $C_5$, $P_5$ and $\bar P_5$ and the \emph{spider and quasi-spider graphs with $P_4$-tidy heads} defined below \cite{Giako}.
We will analyze the behaviour of the Grundy total dominating sequences for these particular classes of modular graphs.

A graph is a \emph{spider graph} if its vertex set can be partitioned into three sets $S$, $C$ and $H$ ($H$ possible empty), where $S$ is a stable set, $C$ is a clique, $|S|=|C|=r \geq 2$, $H$ is completely joined to $C$, and no vertex of $H$ is adjacent to a vertex in $S$. Moreover, if $S=\{s_1, \ldots, s_r\}$ and $C= \{c_1, \ldots, c_r\}$ one of the following conditions must holds:
\begin{enumerate}
\item \emph{thin spider}: $s_i$ is adjacent to $c_j$ if and only if $i = j$.
\item \emph{thick spider}: $s_i$ is adjacent to $c_j$ if and only if $i \neq j$.
\end{enumerate}
The size of $C$ (and $S$) is called the \emph{weight} of $G$ and the set $H$ in the partition is called \emph{the head} of the spider. A spider graph $G$ with the partition $S,C,H$ will be denoted $G=(S,C,H)$.

Notice that if $r=2$ thin and thick spider graphs are both $P_4$. In what follows we consider thick spider graphs with $r\geq 3$.

Now, if $G=(S,C,H)$ is a thin (resp. thick) spider graph, the graph obtained by replacing one vertex $v\in S\cup C$ by $\overline K_2$ or $K_2$ is called \emph{thin (resp. thick) quasi-spider graph}. without loss of generality we assume that the vertex replaced in $S$ is $s_r$ and the vertex replaced in $C$ is $c_r$. Note that the replacement by $\overline K_2$ ($K_2$) is equivalent to add a false (true) twin vertex $s'_r$ or $c'_r$ of $s_r$ or $c_r$, respectively.
We denote by $(S\hookleftarrow W,C,H)$ and $(S,C\hookleftarrow W,H)$ the quasi-spider graph obtaining from a spider graph with partition $(S,C,H)$ by replacing one vertex in $S$ by $W\in\{\overline K_2,K_2\}$ or one vertex in $C$ by $W\in\{\overline K_2,K_2\}$, respectively. 
The weight of a quasi-spider graph is the weight of the original spider graph. 

It is known that the partition for spider and quasi-spider graphs is unique and its recognition as well as its partition can be performed in linear time (see \cite{Giako}).

From Proposition \ref{cor:twins} we infer that $\grt (S\hookleftarrow \overline K_2,C,H)=\grt (S,C\hookleftarrow \overline K_2,H)=\grt (S,C,H)$. Thus it only remains to compute the Grundy total domination number for spider graphs and quasi-spider graphs of the type $(S\hookleftarrow K_2,C,H)$ and $(S,C\hookleftarrow K_2,H)$.

\begin{prop}  \label{thinspider}  
Let $G=(S,C,H)$ be a thin spider graph of weight $r\geq 2$ and $T$ a Grundy total dominating sequence of $G[H]$. Then: 
\begin{enumerate}
\item $T'=(s_1,\dots,s_r)\oplus T\oplus (c_1,\dots,c_r)$ is a Grundy total dominating sequence of $G$ and $(S,C\hookleftarrow K_2,H)$. 
Besides, $$\grt (G)=\grt (S,C\hookleftarrow K_2,H)=\grt (G[H])+2r.$$

\item If $G[H]$ has an isolated vertex $v$, then $T'=(s_1,\dots,s_{r-1})\oplus (v)\oplus T\oplus (c_1,\dots,c_{r-1})\oplus (s_r)\oplus (c_r)$ is a Grundy total dominating sequence of $(S\hookleftarrow K_2,C,H)$. Otherwise, $T'=(s_1,\dots,s_r)\oplus T\oplus (c_1,\dots,c_r)$ is a Grundy total dominating sequence of $(S\hookleftarrow K_2,C,H)$. Besides,
$$\grt (S\hookleftarrow K_2,C,H)=\grt (G)+2 \eta(G[H]).$$
\end{enumerate}
\end{prop}

\proof
\begin{enumerate}
\item Observe that $T'$ is a legal sequence of $G$ and $(S,C\hookleftarrow K_2,H)$.
 
The equality $\grt (G)=\grt (G[H])+2r$ follows immediately from Proposition~\ref{cor:pendant}.

Now consider the thin quasi-spider graph $(S,C\hookleftarrow K_2,H)$. From Proposition~\ref{cor:pendant}, it is enough to analyze the graph $G'$ obtained by deletion of the $r-1$ pendant vertices $s_1,\ldots, s_{r-1}$ and their neighbours. It is easy to see that $G'$ is the join of $K_2$ and the graph obtained by the disjoint union of $G[H]$ and the trivial graph having $s_r$ as its vertex. By Proposition \ref{join},
$\grt (G')=max \{2, \grt (G[H])+2\}= \grt (G[H])+2$. Then,
$\grt (S,C\hookleftarrow K_2,H)=\grt (G[H])+2r= \grt (G)$.

\item
Observe that $T'$ is a legal sequence of $(S\hookleftarrow K_2,C,H)$.
 
Following the same reasoning that in the previous item, in this case $G'$ is the join of the graph $G[H]\oplus K_2$ and the trivial graph having $c_r$ as its vertex. By Proposition \ref{join},
$\grt (G')=\grt (G[H]\oplus K_2)+2 \eta(G[H]\oplus K_2)= \grt (G[H])+2 \eta(G[H])+2$. Then, $\grt (S\hookleftarrow K_2,C,H)= \grt (G[H])+2 \eta(G[H])+2 r=\grt (G)+2 \eta(G[H])$.
\end{enumerate}
\qed

Considering thick spider and quasi-spider graphs we have the following result.

\begin{prop}  \label{thickspider}  
Let $G=(S,C,H)$ be a thick spider graph of weight $r\ge 3$ and $Z$ a Grundy total dominating sequence of $G[H]$. Then, 
\begin{enumerate}
\item $Z'=(s_1,s_2)\oplus Z\oplus (c_1,c_2)$ is a Grundy total dominating sequence of $G$ and $(S,C\hookleftarrow K_2,H)$. Besides,
$$\grt(G)=\grt (S,C\hookleftarrow K_2,H)= 4+\grt (G[H]).$$
\item	$Z'=(s_1,s_2,s_r,s'_r)\oplus Z \oplus (c_1,c_2)$ is a Grundy total dominating sequence of $(S\hookleftarrow K_2,C,H)$. Besides,
$$\grt (S\hookleftarrow K_2,C,H)=\grt (G)+2=6+\grt (G[H]).$$
\end{enumerate}
\end{prop}

\proof
\begin{enumerate}
	\item 

We obtain the lower bound from the fact that the sequence $Z'$ is a Grundy total dominating sequence of $G$ and $(S,C\hookleftarrow K_2,H)$.

For the upper bound, we first consider the thick spider graph $G$. 

Let $T$ be a legal sequence of $G$.
Since for all $i\neq j$, $N(c_i) \cup N(c_j)=V(G)$, $|\widehat T\cap C|\leq2$. Let $I$ be the set of isolated vertices of $G[H]$. It is not hard to see that $|\widehat T\cap (S\cup I)|\leq 2$. Finally, since $|\widehat T\cap (H\setminus I)|\leq\grt (G[H])$, $\grt(G)\leq 4+\grt (G[H])$.

Now, consider $T$ a legal sequence of the thick quasi-spider graph $(S,C\hookleftarrow K_2,H)$. As before, $|\widehat T\cap (S\cup I)|\leq 2$ and $|\widehat T\cap (H\setminus I)|\leq\grt (G[H])$. If $|\widehat T\cap \{c_1,c_2,\dots,c_r,c'_r\}|\leq2$ then $\grt (S,C\hookleftarrow K_2,H)\leq 4+\grt (G[H])$.

If $|\widehat T\cap \{c_1,c_2,\dots,c_r,c'_r\}|>2$ then $c_r$ and $c'_r$ belong to $\widehat T$ and $|\widehat T\cap \{c_1,c_2,\dots,c_r,c'_r\}|=3$. Let $i\in\{1,\dots,r-1\}$ such that $\{c_r,c'_r,c_i\}=\widehat T\cap \{c_1,c_2,\dots,c_r,c'_r\}$. Note that $c_r$ and $c'_r$ must appear before $c_i$ in $T$. without loss of generality we assume that $c'_r$ appears after $c_r$ in $T$. Therefore, since $\{c_1,c_2,\dots,c_r,c'_r\}\subseteq N(c_r)\cup N(c'_r)$, there is no vertex of $I\cup S$ in $T$ after $c'_r$. Besides, since $c'_r$ only footprints $c_r$, the only possible vertex of $I\cup S$ in $T$ before $c'_r$ is $s_r$. Then, $|\widehat T\cap (S\cup I)|\leq 1$ and thus $\grt (S,C\hookleftarrow K_2,H)\leq 4+\grt (G[H])$.

\item
Let us consider the thick quasi-spider graph $(S\hookleftarrow K_2,C,H)$. Clearly, by deleting the true twin vertex $s'_r$ we obtain $G$. Then, by Lemma \ref{grundysubgraph}, we have  $\grt (S\hookleftarrow K_2,C,H)\leq \grt (G)+2$. To complete the proof we only need to observe that $Z'$ 
is a Grundy total dominating sequence of $(S\hookleftarrow K_2,C,H)$.
\end{enumerate}
\qed

As we have already mentioned, non-trivial modular $P_4$-tidy graphs are spider and quasi-spider graphs with $P_4$-tidy heads and the graphs $C_5$, $P_5$ and $\bar P_5$. Therefore, 
combining Propositions \ref{thinspider} and \ref{thickspider}, and the aproach given in Section~\ref{sec:modular}, we derive the following theorem.

\begin{thm}
A Grundy total dominating sequence can be obtained in linear time for $P_4$-tidy graphs.
\end{thm}

\section{NP-completeness in split graphs}
\label{sec:NP}

The following problem was studied in~\cite{bgm-2014}.

{\sc }

\begin{center}
\fbox{\parbox{0.87\linewidth}{\noindent
{\sc Grundy Domination Number Problem}\\[.8ex]
\begin{tabular*}{.96\textwidth}{rl}
{\em Input:} & $G=(V,E)\;,\; k\in\mathbb{Z}^+$.\\
{\em Question:} & \text{Is there a Grundy dominating sequence of} $G$ \text{of length at least} $k$?
\end{tabular*}
}}
\end{center}

\medskip

It was shown that {\sc Grundy Domination Number Problem} is NP-complete, even when restricted to chordal graphs. On the other hand, the following problem

\begin{center}
\fbox{\parbox{0.93\linewidth}{\noindent
{\sc Grundy Total Domination Number Problem}\\[.8ex]
\begin{tabular*}{.96\textwidth}{rl}
{\em Input:} & $G=(V,E)\;,\; k\in\mathbb{Z}^+$.\\
{\em Question:} & \text{Is there a Grundy total dominating sequence of} $G$ \text{of length at least} $k$?
\end{tabular*}
}}
\end{center}

\medskip
\noindent was shown to be NP-complete even when restricted to bipartite graphs \cite{bhr-2016}. In the proof a translation from a certain covering problem in hypergraphs was used.
As the next theorem shows, {\sc Grundy Total Domination Number Problem} remains NP-complete, when restricted to split graphs, which contrasts the result from~\cite{bgm-2014}, showing that {\sc Grundy Domination Number Problem} is polynomial time solvable in split graphs. (Recall that a graph $G$ is a {\em split graph}, if its vertex set can be partitioned into two subsets, one of which induces a clique and the other is a stable set.)

\begin{thm}
{\sc Grundy Total Domination Number} is NP-complete, even when restricted to split graphs.
\end{thm}
\proof
It is clear that the problem is in NP. 

Given a graph $G=(V,E)$ with no isolated vertices, we construct the split graph $G'=(V_1 \cup V_2, E')$ as follows: $V_1=\{v^1: v\in V\}$ is a stable set, $V_2=\{v^2: v\in V\}$ induces a clique and $N_{G'}(v^1_i)=\{v^2\in V_2: v\in N_G(v_i)\}$. 

We will prove that $\gamma_{gr}^t (G')= 2 \gamma_{gr}^t (G)$, which by the NP-completeness result on {\sc Grundy Total Domination Number Problem} from~\cite{bgm-2014} readily implies that the problem is NP-complete even when restricted to split graphs.

If $(v_1,\ldots, v_k)$ is a Grundy total dominating sequence of $G$, then  $(v^1_1,\ldots, v^1_k, v^2_1,\ldots, v^2_k)$ is a legal 
sequence of $G'$. Hence $\gamma_{gr}^t (G')\geq 2 \gamma_{gr}^t (G)$.

Now, let  $S=(w_1,\ldots, w_k)$ be a maximal legal sequence of  $G'$. Let  $t=\min\{j: w_j\in V_2\}$ and  $r=\min\{j: j>t, w_j\in V_2\}$. It is clear that, for all $j\geq r$, $w_j\in V_2$. We can thus write $$S=(v_1^1,\ldots,v_{t-1}^1,v_{1'}^2,v_t^{1},\ldots,v_{r-2}^1,v_{2'}^2,\ldots,v_{(k-r+2)'}^2),$$ 
where $w_t=v_{1'}^2$, $w_r=v_{2'}^2$ appear in $S$ as the first and the second vertex from $\widehat{S}\cap V_2$. (Note that we allow $t-1=r-2$, in which case all vertices from $V_2$ appear in $S$ after all vertices from $V_1$.)

Consider the subsequence $(v_1^1,\ldots,v_{r-2}^1)$ of the vertices in $S$ that are taken from $\widehat S\cap V_1$. Clearly, the corresponding sequence $(v_1,\ldots,v_{r-2})$ in $G$ is a legal sequence in $G$ (not necessarily total dominating sequence). Thus $|\widehat S \cap V_1| \leq \gamma_{gr}^t (G)$.  In addition, the sequence $(v_{2'},\ldots,v_{(k-r+2)'})$ in $G$ that corresponds to the sequence $(v_{2'}^2,\ldots,v_{(k-r+2)'}^2)$ of vertices from $S$ that are taken from $\widehat {S}\cap V_2\setminus \{w_t\}$, is a legal sequence of $G$ (not necessarily total dominating sequence). Hence $|\widehat S\cap V_2|\leq \gamma_{gr}^t (G)+1$.

Now, if $N_G(v_{2'})\setminus N_G(v_{1'})\neq \emptyset$, then the sequence $(v_{1'},\ldots,v_{(k-r+2)'})$ in $G$, which corresponds to the sequence of vertices from $\widehat{S} \cap V_2$, is also a legal sequence in $G$, hence $k-r+2\le \gamma_{gr}^t (G)$; or, in other words, $|\widehat S\cap V_2|\leq \gamma_{gr}^t (G)$. Thus, we infer that $|\widehat S|=|\widehat S\cap V_1|+ |\widehat S\cap V_2|\leq 2 \gamma_{gr}^t (G)$. Finally, suppose that $N_G(v_{2'})\setminus N_G(v_{1'})=\emptyset$. This implies that $v_{2'}^2$ footprints only $v_{1'}^2$. But then the sequence $(v_1^1,\ldots,v_{r-2}^1)$ does not footprint $v_{1'}^2$, which implies that $(v_1,\ldots,v_{r-2})$ is not a total dominating sequence of $G$, and so $r-2=|\widehat S\cap V_1|\leq \gamma_{gr}^t (G)-1$. We again infer $|\widehat S|=|\widehat S\cap V_1|+ |\widehat S\cap V_2|\leq 2 \gamma_{gr}^t (G)$, which completes the proof. 
\qed

\section*{Acknowledgement}
This work was in part supported by MINCYT-MHEST SLO 1409, the Slovenian Research Agency (ARRS) under the grants P1-0297, N1-0043 and J1-7110, and the grants PICT-ANPCyT 0586, PIP-CONICET 277, Argentina.


\end{document}